\documentclass{article}
\usepackage{amsfonts}
\usepackage{graphicx}
\usepackage{epstopdf}
\usepackage{amssymb,amsmath,amscd}
\usepackage{mathrsfs}
\usepackage{amsfonts}
\date{}
\newcommand\be{\begin{array}}
\newcommand\en{\end{array}}
\newcommand\di{\displaystyle}

\newcommand\ga{\gamma}
\newcommand\ep{\emptyset}

\newcommand\si{\sigma}

\newcommand\na{\nabla}
\newcommand\vr{\varphi}

\newcommand\Df{D^+(\frac{3}{4}\mu_0)}
\newcommand\Dg{D^-(\frac{3}{4}\mu_0)}
\newcommand\pa{\partial}

\newcommand\va{\varepsilon}
\newcommand\la{\lambda}
\newcommand\Om{\Omega}
\newcommand\wi{\widetilde}
\newcommand\de{\delta}
\newcommand\De{\Delta}
\newcommand\Da{D^+(\mu_0)}
\newcommand\Db{D^-(\mu_0)}
\newcommand\Dc{D^+(\mu_0)\cap D^-(\mu_0)}

\newcommand\gl{\geqslant}
\newcommand\ls{\leqslant}

\newcommand\Br{B(r, \bar \va, \de)}

\newcommand\ri{\rightarrow}
\newcommand\ti{\times}

\newcommand\bk{\backslash}
\newcommand\al{\alpha}
\newcommand\lge{\langle}
\newcommand\rge{\rangle}

\newcommand\no{\eqno}
\newcommand\mb{\mathbb}

\begin{document}
\title{\bf     Sign Changing Critical Points  for Locally Lipschitz Functionals
\thanks{This paper is supported by NSFC 11871250; Qin Baoxia is the corresponding author.}
}
\author{
Xian  Xu$^{1}$,  Baoxia Qin$^2$ } \maketitle

\vspace{2mm} \noindent \begin{center}\ $^1$Department of Mathematics,
Jiangsu Normal University, Xuzhou,
\\Jiangsu,
221116,  P. R. China

$^2$School of Mathematics,
Qilu Normal  University, Jinan,
\\Shandong,
250013,  P. R. China\\

\end{center}
\begin{center}
\begin{minipage}{5in}
{\small {\bf Abstract}\quad   In this paper,  some existence results for sign-changing critical points of locally Lipschitz functionals in real Banach space are obtained by the method combining the invariant sets of descending flow method with a quantitative deformation. First we assume the locally Lipschitz functionals   to be  outwardly directed on the the boundary of some closed convex sets of the real Banach space. By using the relation between  the critical points on the Banach space and those of the closed convex sets, we construct a  quantitative deformation lemma, and then we obtain some linking type of critical points theorems. These theoretical results can be applied to the study of the existence of sign-changing solutions for differential inclusion problems.
In contrast with the related results in the literatures, the main results of this paper relax the requirement that the functional being of $C^1$ continuous  to locally Lipschitz.
}
\end{minipage}
\end{center}

\section{Introduction}

\ \ \ \ \    The main purpose of this paper is to obtain some  existence results for sign-changing  critical points of  Locally Lipschitz functionals  in Banach spaces.  Chang [1], in ordered  to study  equations with discontinuities, developed an extension of the classical smooth critical point theory, to non-smooth locally Lipschitz functionals.   The theory of Chang was based on the sub-differential of  locally Lipschitz functionals due to Clarke  [2]. Using this sub-differential, Chang proposed a generalization of the deformation lemma and obtained a  Mountain Pass Theorem for Locally Lipschitz functionals. Subsequently, many other types critical points theorems for locally Lipschitz functionals have been obtained and been applied to  the studying of the existence of solutions for various of  elliptic equations boundary value with discontinuities; see [3-15] and the references therein.

 In the past 20 years or so, people have extensively studied the existence of  sign-changing solutions for elliptic boundary value problems; See the literatures [19, 29, 30]. At the same time, some authors have also established the theoretical results of the existence of sign-changing critical points for smooth functionals in abstract space; see [20,25,26,27,32,33]. These theoretical results can be applied to the study of sign-changing solutions for elliptic boundary value problems.  Let's first recall some of these results.
Li   and Wang  [20] presented several variants of Ljusternik-Schnirelman type theorems in partially ordered Hilbert spaces, which assert the locations of the critical
points constructed by the minimax method in terms of the order structures.  The method to show their main results is combining the invariant sets of descending flow method with a quantitative deformation. Li and Wang's method have been applied to various of differential equations boundary value problems to obtained the multiplicity of sign-changing solutions by many authors; see [20, 25-28] and the references therein. T. Bartsch and T. Weth [33] established an abstract critical theory in partially ordered Hilbert spaces by virtue of critical groups and give some additional properties for sign-changing critical points.  By using the linking notation introduced by Schechter and
Tintarev and a type of Schauder invariance condition on a neighborhood of the cone and negative cone in the Hilbert space, Schechter and Zou [27] proved some   critical points theorem of linking  type ensuring the existence for sign-changing critical points for a $C^1$ functional on the Hilbert space. In order to prove that the neighborhoods of the cone and negative cone are invariant set of descending flow, Schechter and Zou [27] used the method of Sun and Liu [17] and Sun [16,18].  The abstract result was used to solve semilinear eigenvalue problems and problems with jumping nonlinearities resonant or not, and  with the Fu$\check{c}i$k spectrum.
Zou [32]  established some parameter-depending linking theorems, which allow to produce a bounded
and sign-changing Palais-Smale sequence. For even functionals, a parameter-depending fountain
theorem is also obtained in [32] which provides infinitely many bounded and sign-changing Palais-Smale sequences.  The abstract results of [32] were applied to the Schr$\ddot{o}$dinger equation with (or without)
critical Sobolev exponents to obtain the existence of  sign-changing solutions. The positive and negative solutions were also gained as by-products.  For further theoretical results on sign-changing critical points, readers can refer to Zou 's monograph [25].

Recently, some authors have studied the existence of sign-changing solutions for elliptic boundary value problems with discontinuous nonlinearities; See [14,15] and the references therein. For instance, Iannizzotto A., M. Salvatore A. and  Motreanu D. [14] investigated the homogeneous Dirichlet problem for a partial differential inclusion involving the $p$-Laplace operator and depending on a parameter $\la>0$.  The existence of three smooth solutions, a smallest positive, a biggest negative, and a nodal one, is obtained for any $\la$ sufficiently large by combining variational methods with truncation techniques.

 As far as we know, no one has theoretically studied the existence of the sign-changing critical points for locally Lipschitz functionals. Now, a natural question is whether we can generalize the theoretical existence results about the sign-changing critical point of smooth functionals to locally Lipschitz functionals. This paper will serve to fulfill this purpose. It is easy to see that the main difficulty faced by this research is how to establish the descending flow that ensuring cone and the negative cone being invariant. In order to overcome this difficulty, we will introduce the conditions of outwardly directed on convex closed sets and the conditions of the Schauder invariance for  set-valued mapping  as in [11]. We will use the relation theorem between the critical points of the whole space and the critical points of the closed convex set in  [11] when the above conditions hold for closed convex sets. Also, we need to introduce a new class of Palais-Smale conditions on closed convex sets. Based on these theorems and concepts, and based on the method of constructing pseudo-gradient fields on closed convex sets in  [6], we establish a quantitative deformation lemma that asserts the neighbourhood of the cone and the negative cone are invariant. On this basis, we obtain the existence result for the sign-changing critical points in the presence of a linking structure, and apply this abstract result to the study of the sign-changing solution of differential inclusion problems.

 The main results of this paper can be thought as a generalization of some main results in Li  and Wang  [20], Zou [25-27] and  Qian [28].  Comparing the related results in Li and Wang  [20], Zou [25-27] and Qian [28], the main results of this paper relax the requirement that the functional being of $C^1$ to locally Lipschitz. Different from the above literature, we use a different method here.  Obviously, the results here can also be applied to study the existence of sign-changing critical points for smooth functionals.

\vskip 0.1in

\section{Preliminaries }

\ \ \ \ \ \ \   First let us recall some theories concerning the  subdifferential theory of locally Lipschitz functionals due to Clarke [2].  Let $X$ be a real Banach space and $X^*$ its topological dual. A functional $\vr :X\ri \mb R$ is said to be locally Lipschitz, if for everey $x\in X$, there exists a neighbourhood $U$ of $x$ and a constant $k>0$ depending on $U$ such that $|\vr (z)-\vr(y)|\ls k\|z-y\|$ for all $z,y\in U$.  For such a  functional we define generalized directional derivative $\vr^0(x; h)$ at $x\in X$ in the direction $h\in X$ by
$$\vr ^0(x;h)=\lim\limits_{x'\ri x}\sup\limits_{\la\downarrow 0^+}\frac{\vr(x'+\la h)-\vr (x')}{\la}.$$
The function $h\mapsto \vr^0(x:h)$ is sublinear, continuous. So by the Hahn-Banach theorem we know that $\vr^0(x:\cdot)$ is the support function of a nonempty, convex and $w^*$-compact set
$$\pa\vr (x)=\{x^*\in X^*:\lge x^*,h\rge\ls \vr^0(x;h)\ \mbox{for all}\ h\in X\}.$$
The set $\pa\vr(x)$ is called the generalized or Clarke subdifferential of $\vr$ at $x$.

If $\vr,\psi:X\mapsto\mb R$ are locally Lipschitz functionals, then $\pa(\vr+\psi)(x)\subset \pa\vr (x)+\pa\psi(x)$, while for any $\la\in \mb R$ we have $\pa(\la\vr)(x)=\la\pa\vr(x)$.  Moreover, if $\vr:X\mapsto \mb R$ is also convex, then this sudifferential coincides with the subdifferential in the sense of convex analysis. If $\vr$ is strictly differentiable, then $\pa\vr(x)=\{\vr'(x)\}$. A point $x\in X$ is a critical point of  $\vr$ if $0\in \pa\vr(x)$.

  S.T. Kyritsi and N. S. Papageorgiou  [6] developed a  critical  point theory for non-smooth locally Lipschitz functionals defined  on closed, convex set  extending this way the work of Struwe. Let $D$  be a nonempty, convex, closed subset of $X$, and $\vr: D\mapsto\mb R$ be a locally Lipschitz  continuous functional.  For $x\in D$ we define
 $$m_D(x)=\inf\limits_{x^*} \sup\limits_{y}\big\{\lge x^*,x-y\rge:y\in D,\|x-y\|<1,x^*\in\pa\vr(x)\big\}.$$

 Evidently, $m_D(x)\gl 0$ for all $x\in D$. This quantity can be viewed as a  measure of the generalized slope of $\vr$ at $x\in D$. If $\vr$ admits an extension $\hat\vr\in C^1(X)$, then  $\pa\vr(x)=\{\vr'(x)\}$ and so we have
 $$m_D(x)=\sup\big\{\lge \vr'(x),x-y\rge:y\in D,\|x-y\|<1\big\},$$
 which is the quantity used by Stuwe [31 ,p.147]. Also if $D=X$, then we have
 $$m_D(x)=m(x)=\inf\{\|x^*\|_{*}: x^*\in \pa\vr(x)\},$$
 which is the quantity used by Chang [1].  According to [6], a point $x\in D$ is called a critical point of $\vr$ on $D$, if $m_D(x)=0$.  It follows from [6] that $m_D: X\mapsto [0,+\infty)$ is a lsc  (lower semi-continuous) function.

 Now we recall some theories concerning convex functionals. Let $\psi :X\mapsto \mb R\cup \{+\infty\}$ be convex, proper and lsc. The function $\psi$ turns out to be  continuous on $D_\psi$, where, as usual
 $$D_\psi:=\{x\in X: \psi(x)<+\infty\}.$$
 If $\pa\psi(x)$ indicates the subdifferential of $\psi$ at the point $x\in X$, namely,
 $$\pa \psi(x):=\{x^*\in X^*:\psi (z)-\psi(x)\gl \lge x^*,z-x\rge,\forall z\in X\},$$ and
 $$D_{\pa\psi}:=\{x\in X:\pa\psi(x)\neq 0\},$$ then $D_{\pa\psi}\subset D_\psi$.  Let $\de_D:X\ri \mb R\cup \{+\infty\}$ be the indicator function of $D$, namely
 $$\de_D(x):=\left\{\be{ll} 0, &\mbox{if}\ x\in D,\\
 +\infty,&\mbox{otherwise}.
 \en\right.
 $$
 Then we have
 $$\pa\de_D(x)=\big\{x^*\in X^*:\lge x^*,z-x\rge\ls 0, \forall z\in D\big\}.$$
The  set $\pa\de_D(x)$ is usually called normal cone to $D$ at $x$. If for  $x\in \pa D$,
$$(-\pa\de_D(x))\cap \pa \vr(x)\subset\{0\},$$
then  we say that $\pa \vr$ turns out to be outwardly directed at $x\in \pa D$. This clearly rewrites as
$$\forall z^*\in \pa\vr(x)\backslash\{0\}\ \mbox{there exists}\ z\in D\ \mbox{fulfilling}\ \lge z^*, z-x\rge<0.\no(2.1)$$
It follows from [11, Lemma 4.1 and 4.2]  we have the following Lemma 2.1 and 2.2.
\vskip 0.1in
 {\bf Lemma 2.1.}  $m_D(x)=0$ if and only if $0\in \pa\vr (x)+\pa\de_D(x)$, namely $\pa\vr (x)\cap (-\pa\de_D(x))\neq\emptyset$.
 \vskip 0.1in

 {\bf Lemma 2.2.} Let $D$ be a closed, convex set of $X$ and $\pa\vr$ be outwardly directed at all points of $\pa D$. If $x_0\in D$ and $m_{D}(x_0)=0$ then $x_0$ is a critical point of $\vr$.
 \vskip 0.1in

 Now let us introduce the Schauder type condition  for set value mappings in a manner as [11]. Let $X$ be  reflexive. As usual, we will identify $X^{**}$ with $X$ while $\mathscr F:X^*\mapsto 2^X$  will denote the duality map, given by
 $$\mathscr F(x^*):=\{x\in X:\lge x^*,x\rge=\|x^*\|^2_{*}=\|x\|^2\}, \forall x^*\in X^*.$$
The set $\mathscr F(x^*)$ turns out to be nonempty, convex, and closed; see,e.g. [5, pp. 311-319]. Define
$$\nabla \vr(x):=\mathscr F(\pa\vr(x)),\ \  x\in X. \no(2.2)$$
Clearly, $\nabla \vr(x)$ depends on the choice of the duality pairing between $X$ and $X^*$ whenever it is compatible with the topology of $X$. If $X$ is a Hilbert space, the duality paring becomes the scalar product and (2.2) gives the usual gradient. Write $I$ for the identity operator on $X$. It follows from [11, Theorem 4.5] we have the following Lemma 2.3.
\vskip 0.1in

 {\bf Lemma 2.3 (Schauder invariance condition).}   Suppose $X$ is reflexive, $D$ is a convex and closed set of $X$, $x_0\in \pa D$ and,
 $$\big(I-\na\vr\big)(x_0)\subset D.\no (2.3)$$
 Then $\pa\vr$ is outwardly directed at $x_0$.
 \vskip 0.1in

 {\bf Remark 2.1.}\ The well known Schauder invariance condition for a $C^1$-functional $\vr$ on a Hilbert space $X$ reds as $(I-\vr')(C)\subset C$; see [16-18,21,22]. It has been extend to  Banach spaces in [22]. The notation of Schauder invariance condition was firstly put forward by Sun Jingxian in [17].  As pointed out in [11], the geometric condition that $\pa\vr$ is outwardly directed at $x_0\in \pa D$ is more general than  (2.3).

 \vskip 0.1in

 {\bf Corollary 2.1.} Suppose $X$ is reflexive, $D$ is a convex and closed set of $X$,   $\pa D=\Gamma_1\cup\Gamma_2$. Moreover,
  $$\big(I-\na\vr\big)(x)\subset D,\ \ \ \  \forall x\in \Gamma_1,$$
 and
 $$\pa\vr (x)\cap\big(-\pa\de_{D}(x)\big)\subset\{0\},\ \ \ \ \forall x\in \Gamma_2. $$
  Then $m_{D}(x_0)=0$ if only if $m(x_0)=0$ for any $x_0\in D$.
 \vskip 0.1in

{\bf Lemma 2.4([34],  Von Neumann).}\   Let $X, Y$ be two Hausdorff topological linear spaces, $C\subset X$, $D\subset Y$ be two convex and compact sets. Let $\psi: X\ti Y\mapsto \mb R$ satisfy:
1)\ $x\mapsto \psi(x, y)$ is upper semi-continuous (usc.) and concave; 2)\ $y\mapsto \psi(x, y)$ is lower semi-continuous (lsc.) and convex.
Then $\psi$ has at least one saddle point $(\bar x, \bar y)\in C\ti D$, that is
$$\psi(x, \bar y)\ls \psi(\bar x, \bar y)\ls \psi(\bar x,  y)\ \mbox{for}\ (x, y)\in C\ti D.$$
\vskip 0.1in

{\bf Definition 2.1.}\  Let $D$ be a nonempty closed subset of $X$, $D_1\subset D$, $r\in\mb R$. We say that $\vr$ satisfies  the non-smooth $C$-condition on $D_1$ with respect to $D$ at level $r$, denoted by $(PS)_{D,D_1}^r$, if every sequence $\{x_n\}\subset D_1$ such that $\vr(x_n)\ri r$ and $(1+\|x_n\|) m_D(x_n)\ri 0$ as $n\ri \infty$, has a convergent subsequence. If $\vr$ satisfies  the non-smooth $C$-condition on $D_1$ with respect to $D$ at every level $r\in \mb R$, we say $\vr$ satisfies  the non-smooth $C$-condition on  $D_1$ with respect to $D$, denoted by $(PS)_{D,D_1}$. We say $\vr$ satisfies  the non-smooth $C$-condition on $D_1$ if $D=X$, denoted by $(PS)_{D_1}$.
\section{Main Results}

\ \ \ \
Let $X$ be a real reflexive Banach space  and $J: X\mapsto \mb R$ be locally Lipschitz.  Let $P$ be a  cone of $X$, that is, $P$ is closed convex set in $X$, $\la x \in P$ for all $x\in P$ and $\la\gl 0$, and $P\cap (-P)=\{0\}$. For each $D\subset X$, denoted by $\mbox{int} D$ the interior of $D$. Let $x_0$ be a critical point of $J$. Then, we call $x_0$  a positive, or negative, or sign-changing critical point of $J$, if $x_0\in P$, or  $x_0\in (-P)$, or $x_0\in X\bk\big(P\cup (-P)\big)$, respectively.  For each $\mu>0$, let
$$D^{\pm}(\mu)=\{ x\in X: \mbox{dist}\ (x, \pm P)\ls \mu\}.$$ Obviously, $D^\pm(\mu)$ are closed convex set in $X$, whose interiors are non-empty. Let us denote by $\pi_P: X\mapsto P$ the orthogonal projection on the cone $P$ so that $T_P(u)=u-\pi(u)$ realizes the distance of $u$ from $P$, that is $\mbox{dist} (u, P)=\|T_P(u)\|$.  Since $X$ is reflexive,  the the distance of $u$ from $P$  is achievable.
\vskip 0.1in

We will employ the following conditions in this paper.

(H$_1)$\ $(I-\na J) \big(\pa D^\pm (\mu_0)\big)\subset D^\pm (\frac{1}{2}\mu_0)$ for some $1>\mu_0>0$;

(H$_2)$\ $J$ satisfies the conditions $(PS)_{X}$, $(PS)_{\Da,D^+(\frac{3}{4}\mu_0)}$, $(PS)_{\Db,D^-(\frac{3}{4}\mu_0)}$, and \\ $(PS)_{\Da\cap \Db}$.

\vskip 0.1in

{\bf Remark 3.1.}\ Obviously, the Schauder invariance condition holds for $D^\pm(\mu_0)$ if the condition (H$_1)$ holds. The condition (H$_1)$ is a variant of the so-called $K$-invariant condition for  the functional being locally Lipschitz.  The $K$-invariant condition for the  the functional beng of $C^1$ class was put forward by Monica Conti, et.al [24]. Some authors have employed  this condition to show the multiplicity results for solutions of various elliptic differential equations in the presence of lower and upper solutions, and some authors have employed  this condition to show the existence of sign-changing solutions; see [25-27].
\vskip 0.1in

{\bf Lemma 3.1.}\ Suppose that (H$_1)$ holds.  Then  the following conclusions hold:

(1)\  If $m_{D^+(\mu_0)}(x_0)=0$ for some $x_0\in D^+(\mu_0)\bk D^-(\mu_0)$, then $m(x_0)=0$;

(2)\  If $m_{D^-(\mu_0)}(x_0)=0$ for some $x_0\in D^-(\mu_0)\bk D^+(\mu_0)$, then $m(x_0)=0$;

(3)\  If $m_{D^+(\mu_0)\cap D^-(\mu_0)}(x_0)=0$ for some $x_0\in D^+(\mu_0)\cap D^-(\mu_0)$, then $m(x_0)=0$.

{\bf Proof.}\  The the conclusions (1) and (2)  can be easily proved by Lemma 2.2. We need only to  show the conclusion (3).

Now we show that $\pa J$  is outwardly directed on $\pa (\Da\cap \Db)$, that is
$$\pa J(x)\cap \big(-\de_{\Da\cap \Db}(x)\big)\subset \{0\},\ \ \forall x\in \pa\big(\Da\cap\Db\big).\no(3.1)$$
Obviously, we have
$$\pa\big(\Da\cap\Db\big)=\big(\pa\Da\cap \Db\big)\cup \big(\pa \Db\cap\Da\big).$$
Let $x_0\in \pa\big(\Da\cap\Db\big)$. Assume  that $x_0\in \pa\Da\cap \Db$. Then we have following two cases:

(a)\  $x_0\in \pa\Da\cap \mbox{int}(\Db)$. We claim that $$\pa \de_{\Dc}(x_0)\subset \pa\de_{\Da}(x_0).$$
In fact, for any $x^*\in \pa\de_{\Dc}(x_0)$, by the definition we have
$$\lge x^*, z-x_0\rge \ls 0, \ \ \forall z\in \Dc.$$
Since $x_0\in \pa\Da\cap \mbox{int}(\Db)$,  there exists $r_0>0$ such that
$$B(x_0, r_0)\cap\Da\subset \Dc.$$ For each $z\in \Da\bk B(x_0, r_0)$,  we have
$$z_1:=x_0+r_0\frac{z-x_0}{\|z-x_0\|}=\Big(1-\frac{r_0}{\|z-x_0\|}\Big) x_0+\frac{r_0}{\|z-x_0\|} z\in \Da,$$
and $\|z_1-x_0\|= r_0$. Thus, $$z_1\in B(x_0, r_0)\cap \Da\subset \Dc.$$ So, we have
$$0\gl \lge x^*, z_1-x_0\rge =\frac{r_0}{\|z- x_0\|}\lge x^*, z-x_0\rge. $$
Consequently, we have $\lge x^*, z-x_0\rge\ls 0$ for all $z\in \Da$. This implies that $x^*\in \pa\de_{\Da} (x_0)$. Hence,
$$\pa\de_{\Da} (x_0)\supset\pa\de_{\Dc} (x_0).\no(3.2)$$

By (H$_1)$,  for any $z^*\in \pa J(x_0)\bk\{0\}$ and $y\in \mathscr F(z^*)$, there exists $z\in \Da$ such that $x_0-y=z$. So,
$$\lge  z^*, z- x_0\rge=-\lge z^*, y\rge=-\|z^*\|^2_{*}<0. $$
Hence,
$$\pa J(x_0)\cap \pa \de_{\Da}(x_0)\subset \{0\}.\no(3.3)$$
It follows from (3.2) and (3.3) that (3.1) holds.

(b)\  $x_0\in \pa\Da\cap \pa\Db$. By (H$_1)$, we have $$x_0-\na J(x_0)\subset \Dc.$$ So, for all $z^*\in \pa J(x_0)\bk\{0\}$ and $y\in \mathscr F(z^*)$, there exists $z\in \Dc$ such that $x_0-y=z$. Then we have
$$\lge z^*, z-x_0\rge=-\lge z^*, y\rge=-\|z^*\|^2_{*}<0.$$
This implies  (3.1) holds.

It follows from Lemma 2.2 that $m(x_0)=0$  if $m_{\Dc}(x_0)=0$. The proof is complete.

\vskip 0.1in

Let $r, r_1, r_2\in \mb R$ with $r_1\ls r_2$. In what follows we use the following notations:
$$\mb K=\{x\in X: m(x)=0\},  \mb K_r=\{x\in X: m(x)=0\ \mbox{and}\ J(x)=r\},$$
$$ (\mb K_r)_\de=\{x\in X:\mbox{dist}\ (x, \mb K_r)<\de\},  J^r=\{ x\in X: J(x)\ls r\},$$
$$J_r=\{ x\in X: J(x)\gl r\}, J_{r_1}^{r_2}=\{x\in X: r_1\ls J(x)\ls r_2\}.$$
\vskip 0.1in

{\bf Lemma 3.2.}\ Suppose that (H$_1)$   and (H$_2)$ hold.  Then, for any $\de>0$, there exist $b, \bar\va>0$ such that
$$(1+\|x\|)m(x)\gl b, \forall x\in B(r, \bar\va, \de)\bk (\Da\cup \Db), \no(3.4)$$
$$(1+\|x\|)m_{\Da}(x)\gl b, \forall x\in B(r, \bar\va, \de)\cap (\Da\bk \Db), \no(3.5)$$
$$(1+\|x\|)m_{\Db}(x)\gl b, \forall x\in B(r, \bar\va, \de)\cap (\Db\bk \Da), \no(3.6)$$
$$(1+\|x\|)m_{\Dc}(x)\gl b, \forall x\in B(r, \bar\va, \de)\cap \Da\cap \Db, \no(3.7)$$
where $\Br=J^{r+\bar\va}_{r-\bar\va}\bk(\mb K_r)_\de$.

{\bf Proof.}\ First we prove that there exist $\bar\va_1, b_1>0$ such that
$$(1+\|x\|)m(x)\gl b_1, \forall x\in B(r, \bar\va_1, \de)\bk \big(\Df\cup \Dg\big).\no(3.8)$$
Arguing by make contradiction that this is not the case. Then there exist positive numbers sequences $\{\bar b_n\}$, $\{\bar \va_n\}$ and $\{x_n\}\subset B(r, \va_n,\de)\bk (\Df\cup \Dg)$ such that $\bar b_n\ri 0$, $\bar \va_n\ri 0$ as $n\ri\infty$, and
$$(1+\|x_n\|) m(x_n)<\bar b_n.$$
By using the condition  (H$_2)$, there  is a subsequence of $\{x_n\}$ and $x_0\in \overline{X\bk(\Df\cup\Dg\cup (\mb K_r)_\de)}$,  up to a subsequence, we assume that $x_n\ri x_0$ as $n\ri\infty$. Obviously, we have $J(x_0)=r$ and $m(x_0)=0$, that is $x_0\in \mb K_r\subset (\mb K_r)_\de$. Since $\{x_n\}\subset B(r, \va_n,\de)\bk (\Df\cup \Dg)$ for $n=1,2,\cdots$,  we have $x_0\not\in (\mb K_r)_\de$, which is a contradiction. Hence, (3.8) holds.

 It follows from (H$_1)$  that $$\mb K_r\cap \big(D^+(\mu_0)\bk(D^-(\mu_0)\cup D^+(\frac{3}{4}\mu_0))\big)=\ep. $$  Thus, by (3.8)  we have
$$(1+\|x\|)m(x)\gl  b_1, \  \forall x\in B(r, \bar\va_1, \de)\cap \big(D^+(\mu_0)\bk(D^-(\mu_0)\cup D^+(\frac{3}{4}\mu_0))\big).$$
For each $x_0\in B(r,  \bar \va_1,\de)\cap \big(D^+(\mu_0)\bk(D^-(\mu_0)\cup D^+(\frac{3}{4}\mu_0))\big)$ and $x^*\in \pa J(x_0)$, by (H$_1)$ we may assume that $\mathscr F x^*=x_0-w$ for some $w\in  D^+(\frac{1}{2}\mu_0)$. Let
$\bar x=x_0+\frac{\mu_0}{8\|w-x_0\|}(w-x_0)$. Then we have $\|w-x_0\|\gl \frac{1}{4}\mu_0$. In fact, if $\|w-x_0\|<\frac{1}{4}\mu_0$, then we have
$$\be{ll}\frac{3\mu_0}{4}&\ls\mbox{dist}\ (x_0,P)=\|T_P x_0\|\ls \|x_0-\pi_P(w)\|\\
&\ls\|x_0-w\|+\|w-\pi_P(w)\|=\|x_0-w\|+\|T_P(w)\|\\
&< \frac{\mu_0}{4}+\frac{\mu_0}{2}<\frac{3\mu_0}{4},\en$$
which is a contradiction. This implies that $\|w-x_0\|\gl \frac{1}{4}\mu_0$.
Consequently,  we have
$$\bar x=\Big(1-\frac{\mu_0}{8\|w-x_0\|}\Big)x_0+\frac{\mu_0}{8\|w-x_0\|}w\in D^+(\mu_0)$$
and $\|x_0-\bar x\|=\frac{1}{8}\mu_0<1$.  So,
$$\sup\limits_{v\in D^+(\mu_0),\|x_0-v\|<1}\lge x^*, x_0-v\rge\gl \lge x^*, x_0-\bar x\rge=\frac{\mu_0}{8\|w-x_0\|}\lge x^*, x_0-w\rge=\frac{1}{8}\mu_0\|x^*\|_{*},$$
and
$$m_{D^+(\mu_0)}(x_0)\gl \frac{1}{8}\mu_0\inf\limits_{x^*\in\pa J(x_0)}\|x^*\|_{*}= \frac{1}{8}\mu_0 m(x_0).$$
Hence,
$$(1+\|x\|)m_{D^+(\mu_0)}(x)\gl \frac{1}{8}\mu_0 b_1, \  x\in B(r, \bar \va^+,\de)\cap \big(D^+(\mu_0)\bk(D^-(\mu_0)\cup D^+(\frac{3}{4}\mu_0))\big).$$

Now we show that there exist $\wi\va_2, \wi b_2>0$ such that
$$(1+\|x\|)m_{\Da}(x)\gl \wi b_2, \forall x\in B(r, \wi\va_2, \de)\cap (\Df\bk \Db).$$
Arguing by make contradiction that this is not the case. Then there  exist positive numbers sequences $\{\bar b_n\}$, $\{\bar \va_n\}$ and sequence $\{x_n\}\subset B(r, \bar\va_n, \de)\cap (\Df\bk \Db)$ such that $\bar b_n\ri 0$, $\bar \va_n\ri 0$ as $n\ri\infty$, and
$$(1+\|x_n\|) m_{\Da}(x_n)< \bar b_n.$$
Since $J$ satisfies the condition $(PS)_{\Da,D^+(\frac{3}{4}\mu_0)}$, then, up to  a subsequence if necessary, we assume that $x_n\ri  x_0$ as $n\ri \infty$, where $m_{\Da}(x_0)=0$ and $J(x_0)=r$. It follows from Lemma 3.1 that $m(x_0)= m_{\Da}(x_0)=0$, so $x_0\in\mb K_r$. Since $\{x_n\}\subset B(r, \wi\va_n, \de)\cap (\Df\bk \Db)$ for each $n=1,2,\cdots$, we have $x_0\not\in (\mb K_r)_\de$, which is a contradiction. Let $b_2=\min\{\wi b_2, \frac{1}{8}\mu_0 b_1 \}$ and $\bar\va_2=\min\{\bar \va_1, \wi\va_2\}$. Then we have
$$(1+\|x\|)m_{\Da}(x)\gl  b_2, \forall x\in B(r, \bar\va_2, \de)\cap \big(\Df\bk \Db\big), \no(3.9)$$

In the same way as above we can prove that, there exist $b_3, b_4,\bar\va_3,\bar\va_4>0$ such that
$$(1+\|x\|)m_{\Db}(x)\gl b_3, \forall x\in B(r, \bar\va_3, \de)\cap (\Db\bk \Da), \no(3.10)$$
and
$$(1+\|x\|)m_{\Dc}(x)\gl b_4, \forall x\in B(r, \bar\va_4, \de)\cap \Da\cap \Db. \no(3.11)$$

Now  let  $b=\min\{ b_1, b_2, b_3, b_4\}$ and $\bar \va=\min\{\bar\va_1, \bar\va_2,\bar\va_3,\bar\va_4\}$. It follows from (3.8)$\sim$(3.11) that (3.4)$\sim$(3.7) hold. The proof is complete.
\vskip0.1in

 In what follows of this section, for brevity, let $S= \Br\bk(\Da\cup \Db)$.
 \vskip 0.1in

 {\bf Lemma 3.3.}\ Suppose that (H$_1)$ and (H$_2)$   hold. Then there exists $v:\Br\ri X$ such that $v$ is locally Lipschitz continuous, $\|v(x)\|\ls 2(1+\|x\|)$ for all $x\in \Br$ and
  $$x- \frac{1}{1+\|x\|}v(x)\in \Da, \forall x\in \Br\cap \Da,$$
   $$x- \frac{1}{1+\|x\|}v(x)\in \Db, \forall x\in \Br\cap \Db,$$
   $$\lge x^*, v(x)\rge\gl \frac{b}{8}, \ \forall x\in\Br, x^*\in\pa J(x).$$
   Moreover, $v$ is odd if $J$ is even.

{\bf Proof.}\ For each $x_0\in S$, take $w_0^*\in \pa J(x_0)$ such that $m(x_0)=\|w_0^*\|_{*}>0$. Then we have $B(0, \|w_0^*\|_{*})\cap \pa J(x_0)=\ep$, where $B(0, \|w_0^*\|)=\{z^*\in X^*: \|z^*\|_{*}<\|w_0^*\|_{*}\}$. So,  by Hahn-Banach separation theorem, we can find $u_1(x_0)\in X$ with $\|u_1(x_0)\|=1$ such that for all $z^*\in B(0, \|w_0^*\|_{*})$ and $y^*\in\pa J(x_0)$,
$$\lge z^*, u_1(x_0)\rge \ls \lge w_0^*, u_1(x_0)\rge\ls \lge y^*, u_1(x_0)\rge.$$
Recall that $$\|w_0^*\|_{*}=\sup\big\{\lge z^*, u_1(x_0)\rge: z^*\in B(0, \|w_0^*\|_{*})\big\},$$
then we have by Lemma 3.2,
$$\lge y^*, u_1(x_0)\rge \gl \|w_0^*\|_{*}>\frac{b}{2(1+\|x_0\|)}.$$
Since the map $x\mapsto\pa J(x)$ is usc from $X$ into $X^*_w$, we may take an open neighborhood $B_1(x_0, r_1(x_0))$ of $x_0$,  such that
$$\lge  y^*, u_1(x_0)\rge>\frac{b}{4(1+\|y\|)},\ \ \forall y^*\in\pa J(y), y\in U_1(x_0), \no(3.12)$$
where $U_1(x_0)=B_1(x_0, r_1(x_0))\cap\Br$.
Since $x_0\in S\subset \Br\bk(\Da\cup\Db)$, we may take $r_1(x_0)>0$ small enough such that $U_1(x_0)\subset S$.

Pick $x_0\in \Br\cap \big(\Da\bk \Db\big)$, we have  by Lemma 3.2,
$$m_{\Da}(x_0)>\frac{b}{1+\|x_0\|}.$$
 Let $C=(\{x_0\}-\Da)\cap \bar B(0,1)$ and $D=\pa\vr(x_0)$, where $\bar B(0, 1)=\{x\in X: \|x\|\ls 1\}$. Let $X_w$ and $X^*_w$ be the spaces $X$ and $X^*$ furnished their weak topology respectively. Since $X$ is reflexive,  $D$ is  compact in $X^*_w$ and $C$ is  compact  in $X_w$. Obviously, both $C$ and $D$ are convex. Let $\psi: C\ti D\mapsto \mb R$ be defined by $\psi(x, y^*)=\lge y^*, x\rge$ for any $(x, y^*)\in C\ti D$. It follows from Lemma 2.4 that $\psi$ has at least one saddle point $(u_+(x_0), x_0^*)\in C\ti D$, that is
$$\min\limits_{x^*\in D}\max\limits_{x\in C}\lge x^*, x\rge=\lge x_0^*, u_+(x_0)\rge =\max\limits_{x\in C}\min\limits_{x^*\in D}\lge x^*, x\rge, \ \forall x^*\in \pa\vr(x_0), x\in (\{x_0\}-\Da)\cap \bar B(0,1).$$ Thus, we have
$$\lge x_0^*, x\rge\ls\lge x_0^*, u_+(x_0)\rge\ls \lge x^*, u_+(x_0)\rge, \ \forall x^*\in \pa\vr(x_0), x\in (\{x_0\}-\Da)\cap \bar B(0,1).$$   Then,  for any $x^*\in\pa\vr(x_0)$, we have
$$\lge x^*, u_+(x_0)\rge\gl\lge x_0^*, u_+(x_0)\rge= m_{\Da}(x_0)>\frac{b}{2(1+\|x_0\|)}.$$
Again by using the fact that $x\mapsto\pa J(x)$ is usc, we know that there exists an open neighborhood $B_+(x_0, r_+(x_0))$ of $x_0$  such that for any $y\in U_+(x_0)$, $y^*\in\pa J(y)$, $$\lge y^*,u_+(x_0)\rge>\frac{b}{4(1+\|y\|)}.\no(3.13)$$
where $U_+(x_0)=B_+(x_0, r_+(x_0))\cap \Br$. Since $\Br\bk \Db$ is open in $\Br$ and $x_0\in \Br\bk\Db$, we may assume that $U_+(x_0)\subset \Br\bk \Db$.

By using the same method as above, for each $x_0\in \Db\bk \Da$, there exist an open neighborhood of $x_0$ in $\Br$, denoted by $U_-(x_0)$, and $u_-(x_0)$ such that $\|u_-(x_0)\|\ls 1$,  $x_0-u_-(x_0)\in\Db$, $U_-(x_0)\subset \Br\bk \Da$ and
$$\lge y^*,u_-(x_0)\rge>\frac{b}{4(1+\|y\|)}\no(3.14)$$
for any $y\in U_-(x_0)$, $y^*\in\pa J(y)$.

Also, for each $x_0\in \Dc$, there exist $u_2(x_0)\in X$ and  an open neighborhood $U_2(x_0)$ of $x_0$ in $\Br$ with $\|u_2(x_0)\|\ls 1$, such that $x_0-u_2(x_0)\in\Dc$, and
$$\lge y^*,u_-(x_0)\rge>\frac{b}{4(1+\|y\|)}\no(3.15)$$
for any $y\in U_2(x_0)$, $y^*\in\pa J(y)$.

Set
$$\mathscr A_1=\{ U_1(x_0): x_0\in \Br\bk(\Da\cup\Db)\},$$
$$\mathscr A_2=\{ U_+(x_0): x_0\in \Br\cap (\Da\bk\Db)\},$$
$$\mathscr A_3=\{ U_-(x_0): x_0\in \Br\cap (\Db\bk\Da)\},$$
$$\mathscr A_4=\{ U_2(x_0): x_0\in \Br\cap \Db\cap\Da\},$$
and $\mathscr A=\mathscr A_1\cup \mathscr A_2\cup \mathscr A_3\cup \mathscr A_4$. Obviously, $\mathscr A$ is an open cover of $\Br$. By paracompactness we can find a locally finite refinement $\mathscr D=\{V_\al:\al\in\Lambda\}$ and a locally Lipschitz partition of unit $\{\ga_\al:\al\in\Lambda\}$ sub-ordinate to it. For each $\al\in\Lambda$ we can find $x_\al\in\Br$ such that $V_\al\subset U_{i(\al)}(x_\al)$ for some $i(\al)\in\{+, -, 1,2\}$  and $U_{i(\al)}(x_\al)\in \mathscr A$.  To this $x_\al\in\Br$ corresponds  the element $w_\al$ such that $\|w_\al\|\ls 1$, and
$$w_\al=\left\{\be{ll}
u_1(x_\al), &\mbox{if}\ V_\al\subset U_1(x_\al);\\
u_+(x_\al), &\mbox{if}\ V_\al\subset U_+(x_\al);\\
u_-(x_\al), &\mbox{if}\ V_\al\subset U_-(x_\al);\\
u_2(x_\al), &\mbox{if}\ V_\al\subset U_2(x_\al).\en\right.
$$

Now, let $v:\Br\ri X$ be defined by
$$v(x)=(1+\|x\|)\sum\limits_{\al\in \Lambda} \ga_\al(x) (w_\al-x_\al+x).\no(3.16)$$

For each $x\in U_i(x_\al) $ with $i\in \{+,-,1,2\}$ and $x^*\in\pa J(x)$, there exists $L_{i,\al}>0$ such that $\|x^*\|_{*}\ls L_{i,\al}$. We may also assume that $B_i(x_\al, r_i(x_\al))$ has a small radium $r_i(x_\al)>0$ such that $(1+\|x\|)(1+\|x_\al\|)^{-1}\ls 2$ for each $x\in U_i(x_\al)$, and $$0< r_i(x_\al)\ls \min\Big\{1,\frac{b}{16(1+\|x_\al\|) L_{i,\al}}\Big\}.$$
Then, we have $$\|v(x)\|\ls(1+\|x\|)\sum\limits_{\al\in \Lambda} \ga_\al(x) (\|w_\al\|+\|x_\al-x\|)\ls 2(1+\|x\|),$$
and
$$\be{ll}(1+\|x\|)&\big| \sum\limits_{\al\in \Lambda} \ga_\al(x) \lge x^*, x-x_\al\rge\big|\\
&\ls \sum\limits_{\al\in\Lambda}\frac{1+\|x\|}{1+\|x_\al\|}(1+\|x_\al\|)\ga_\al(x)\|x^*\|_{*}\|x- x_\al\|\\
&\ls 2 \sum\limits_{\al\in\Lambda} L_{i,\al}\|x- x_\al\|(1+\|x_\al\|)\ga_\al(x)\\
&<\frac{b}{8}\en
$$
for each $x\in \Br$.  On the other hand, by (3.12)$\sim$(3.16)  we have $$\sum\limits_{\al\in\Lambda}\ga_\al(x)(1+\|x\|)\lge x^*, w_\al\rge\gl \frac{b}{4}.$$
So,  we have
$$\be{ll} \lge x^*, v(x)\rge&= \sum\limits_{\al\in\Lambda}\ga_\al(x)(1+\|x\|)\lge x^*, w_\al\rge\\
&+\sum\limits_{\al\in\Lambda}\ga_\al(x)(1+\|x\|)\lge x^*, x-x_\al\rge\gl \frac{b}{8}.\en$$

Next we show that, $x-\frac{1}{1+\|x\|}v(x)\in \Da$ for any $x\in \Da\bk\Db$. In fact,
$$x-\frac{1}{1+\|x\|} v(x)=\sum\limits_{\al\in\Lambda}\ga_\al(x)(x_\al- w_\al).$$
If $\ga_\al(x)\neq 0$ for some $\al\in \Lambda$, then $x\in V_\al\subset U_{i(\al)}(x_\al)$ for some $i(\al)\in\{+,-,1,2\}$. As $x\in \Da\bk\Db$, by the definition of $\mathscr A_1, \mathscr A_2, \mathscr A_3, \mathscr A_4$,  we must have $i(\al)\in\{+,2\}$. So, $x_\al- w_\al\in\Da$. By (3.16) and using the  convexity of $\Da$, we have $$x-\frac{1}{1+\|x\|} v(x)\in\Da.$$

 Similarly we can show that $x-\frac{1}{1+\|x\|} v(x)\in\Db$ for each $x\in \Db\bk\Da$, and $x-\frac{1}{1+\|x\|} v(x)\in\Dc$ for each $x\in \Dc$.

 At last we show that $v$ is odd whenever $J$ is even. Obviously, $\Br$ is a symmetric closed set  with respect to 0. Assume without loss of generality that $\Br\cap\Da\neq\ep$. Then we have $\Br\cap \Db\neq \ep$. Let $v:\Br\ri X$ be defined by (3.16) and let
$$v_0(x)=\frac{1}{2}\big[ v(x)- v(-x)\big], \ \forall x\in\Br.$$
Obviously, $v_0:\Br\mapsto X$ is locally Lipschitz, and  odd map, such that $\|v_0(x)\|\ls 2(1+\|x\|)$.
For each $x\in \Br$ and $x^*\in\pa J(x)$, we have
$$\lge x^*, v_0(x)\rge=\frac{1}{2}\big(\lge x^*, v(x)\rge+\lge -x^*, v(-x)\rge\big)\gl\frac{b}{8}.$$
This implies that $v_0$ is a generalized normalized pseudo-gradient vector field of the
 locally Lipschitz functional $J$.

 By using the symmetric property  of $\Da=-\Db$, we have $-x\in \Db\bk\Da$ for each $x\in \Da\bk\Db$. So, we have
 $$x-\frac{1}{1+\|x\|}v(x)\in\Da\ \mbox{and}\  -x-\frac{1}{1+\|x\|}v(-x)\in\Db$$
 for $x\in\Da\bk\Db$. Therefore, for $x\in\Da\bk\Db$, we have$$-\Big(-x-\frac{1}{1+\|x\|}v(-x)\Big)\in\Da,$$ and so
 $$ x-\frac{1}{2(1+\|x\|)}\big(v(x)-v(-x)\big)=\frac{1}{2}\Big(x-\frac{1}{1+\|x\|}\Big)+ \frac{1}{2}\Big(-\big(-x-\frac{1}{1+\|x\|}v(-x)\big)\Big)\in\Da. $$

 Similarly, we have
 $$x-\frac{1}{2(1+\|x\|)}\big(v(x)-v(-x)\big)\in \Db
 $$
for each $x\in\Db\bk\Da$, and $$x-\frac{1}{2(1+\|x\|)}\big(v(x)-v(-x)\big)\in \Dc$$ for each $x\in \Dc$.  The proof is complete.
\vskip 0.1in

The next theorem (deformation theorem) is the key tool for our main results.
\vskip 0.1in

Obviously, if the condition (H$_1)$ holds then $\mb K_r\cap\pa\big(\Da\cup\Db\big)=\ep$ for any $r\in\mb R$.  Let $\mb K_r^1=\mb K_r\cap\big(\Da\cup \Db\big)$ and $\mb K_r^2= S\cap \mb K_r$. Assume (H$_2)$ holds. Then, for each $r\in\mb R$, $\mb K_r^1$ and $\mb K_r^2$ are compact sets, and so we may take $\de>0$ small such that $(\mb K_r^1)_{3\de}\cap(\mb K_r^2)_{3\de}=\ep$.
\vskip 0.1in

{\bf Lemma 3.4.}\ Suppose that (H$_1)$ and (H$_2)$ hold.  Then for each $\va_0>0$, there exists $0<\bar \va<\va_0$, such that for each $0<\va<\bar\va$, compact set $A\subset J^{r+ \va}$, there exists $\eta\in C([0,1]\ti X, X)$ with  the  following properties:

1)\ $\eta(t, u)=u$ for $t=0$ or $u\not\in J^{-1}(r-\va_0, r+\va_0])\bk(\mb K_r^2)_\de$;

2)\ $\eta(1, A\bk(\mb K_r^2)_{3\de})\subset J^{r-\va}\cup \mbox{int}(\Da)\cup \mbox{int}(\Db)$; and $\eta(1, A)\subset J^{r-\va}\cup \mbox{int}(\Da) \cup\mbox{int}(\Db)$ if $\mb K_r^2=\ep$;

3)\ $\eta(t,\cdot)$ is a homeomorphism of $X$ for $t\in [0,1]$;

4)\ $J(\eta(\cdot, u))$ is non-increasing  for any $u\in X$;

5)\ $\eta(t,\Da)\subset \Da$, $\eta(t, \Db)\subset \Db$ and $\eta(t, \Dc)\subset \Dc$ for $t\in [0,1]$; and

6)\ if $J$ is even, $\eta$ is odd in $u\in X$.

{\bf Proof.}\   Since $\mb K_r$ is a compact set, and so $(\mb K_r)_{3\de}$ is a bounded set. Assume that $\|u\|\ls M_0$ for all $u\in (\mb K_r)_{3\de}$ and some $M_0>0$. It follows from  Lemma 3.2 that there exists $\bar\va, b>0$ such that (3.4)$\sim$(3.7) hold. By Lemma 3.3,  there exists a locally Lipschitz mapping $v:\Br\ri X$ such that $\|v(x)\|\ls2( 1+\|x\|)$, and
$$\lge x^*, v(x)\rge \gl\frac{b}{8},\ \  \forall x\in \Br, x^*\in \pa J(x), \no(3.17)$$
$$x-\frac{1}{1+\|x\|}v(x)\in \Da,\ \  \forall x\in \Da, \no(3.18)$$
$$x-\frac{1}{1+\|x\|}v(x)\in \Db, \ \ \forall x\in \Db, \no(3.19)$$
$$x-\frac{1}{1+\|x\|}v(x)\in \Dc,\ \  \forall x\in \Dc. \no(3.20)$$

Without loss of generality, we may assume that $$0<\bar\va<\min\big\{ \va_0,\frac{b\de}{32(1+M_0)}\big\}.$$ For each $0<\va<\bar\va$,  we set
$$D_1=\{x\in X: r-\va\ls J(x)\ls r+\va\}, $$
$$D_2=\{x\in X: r-\bar\va\ls J(x)\ls r+\bar\va\}, $$
and $$\vr(x)=\frac{\mbox{dist}\ (x, X\bk D_2)}{\mbox{dist}\ (x, X\bk D_2)+\mbox{dist}\ (x,  D_1)}.$$
Obviously, $\vr:X\mapsto \mb R$ is locally Lipschitz continuous on $X$. Let
$$\psi(x)=\frac{\mbox{dist}\ (x, (\mb K_r)_\de)}{\mbox{dist}\ (x, (\mb K_r)_\de))+\mbox{dist}\ (x,  X\bk (\mb K_r)_{2\de})}.$$
Also,  $\psi:X\mapsto \mb R$ is locally Lipschitz continuous on $X$.

Consider the following initial value problem (IVP) in $X$
$$\left\{\be{l}\frac{du}{dt}=-V(u),\\
u(0)=u_0,\en\right.\no(3.21)$$
where $V(u)=\vr(u)\psi(u) v(u)$ for $u\in X$. Since $V: X\mapsto X$ is locally Lipschitz continuous on $X$,  by the theory of ordinary differential equations in Banach spaces, (3.21) has
a unique solution in $X$, denoted by $\si(t, u_0)$, with right maximal interval
of existence $[0, T(u_0))$.

Now we show that $T(u_0)=+\infty$ for each $u_0\in X$. Arguing by make contradiction that $T(u_0)<+\infty$. By (3.21) we have
$$\|\si(t, u_0)-u_0\|\ls\int^t_0\|V(\si(s, u_0))\|ds\ls 2\int^t_0\big(1+\|\si(s, u_0)\|\big)ds.$$
So, we have
$$\be{ll}\frac{1}{2}\|\si(t, u_0)-u_0\|&\ls\di\int^t_0(1+\|\si(s, u_0)\|)ds\\
&\ls\di\int^t_0\|\si(s, u_0)-u_0\|ds+(1+\|u_0\|) t,\en$$
By the well known Gronwall's inequality, we have
$$\be{ll}\frac{1}{2}\|\si(t, u_0)-u_0\|&\ls\di\int^t_0(1+\|u_0\|)e^{t-s}ds+(1+\|u_0\|) t\\
&\ls (1+\|u_0\|)(e^{t}-1)+(1+\|u_0\|) t\\
&\ls (1+\|u_0\|)(t+e^{t}-1)\\
&\ls (1+\|u_0\|)(T(u_0)+e^{T(u_0)}).\en$$
So, we have
$$\|\si(t, u_0)\|\ls 2(1+\|u_0\|)(T(u_0)+e^{T(u_0)})+2\|u_0\|=: M_1.$$
Take $\{t_n\}\subset [0, T(u_0))$ such that $t_n\ri T^-(u_0)$ and  for $n=1,2,\cdots$, $$|t_n-t_{n-1}|<\frac{1}{2\cdot 2^n(1+M_1)}.$$
Then we have
$$\be{ll} \|\si(t_n, u_0)-\si(t_{n-1}, u_0)\|&\ls\di\int_{t_{n-1}}^{t_n} \|V(\si(s, u_0)\|ds\\
&\ls 2 \di\int_{t_{n-1}}^{t_n} \big(1+\|\si(s, u_0)\|\big)ds\\
&\ls 2(1+M_1)(t_n-t_{n-1})<\frac{1}{2^n}.\en$$
This implies that $\{\si(t_n, u_0)\}$ is a Cauchy sequence. Thus, there exists $\bar u\in X$ such that $\si(t_n, u_0)\ri \bar u$ as $t\ri T^-(u_0)$.

Now we consider the initial value problem
$$\left\{\be{l}\frac{du}{dt}=-V(u),\\
u(0)=\si(T(u_0), u_0),\en\right.\no(3.22)$$
Then (3.22) has a unique solution on $[0, \bar\de)$ for some $\bar\de>0$ since $V$ is locally Lipschitz on $X$, and so (3.21) has a unique solution on $[0, T(u_0)+\bar\de)$, which is a contradiction. Thus, we have $T(u_0)=+\infty$.

For each $u_0\in X$, denote by $o(u_0)$ the orbit of (3.21) starting from $u_0$, that is
$$o(u_0)=\{\si(t, u_0): t\in [0,+\infty)\}.$$

Now  we show that $\Da$ and $\Db$ are flow invariant sets of (3.21). We only show that $\Da$ is a flow invariant set of (3.21). For a given $u_0\in \pa\Da$, we have the following three cases:

(1) $u_0\in J^{-1}([r-\va, r+\va])\cap \pa \Da$.

In this case, we have
$\vr(u_0)=\psi(u_0)=1$. Thus, for $\la>0$ small enough, by (3.18) and (3.20) we have
$$\be{ll}u_0+\la(-V(u_0))&=\la (1+\|u_0\|)\big(u_0-\frac{1}{1+\|u_0\|}v(u_0)\big)\\
&+(1-\la(1+\|u_0\|))u_0\in\Da.\en$$

(2)\  $u_0\not\in J^{-1}([ r-\bar\va, r+\bar\va])\cap\pa \Da$.

In this case, $V(u_0)=0$. Obviously, for $\la>0$ small enough, by (3.18) and (3.20) we have $u_0+\la(-V(u_0))= u_0\in \Da$.

(3)\ $u_0\in \big(J^{-1}([ r-\bar\va, r+\bar\va])\bk J^{-1}([ r-\va, r+\va])\big)\cap \pa\Da$.

 In this case, $\psi(u_0)=1$ and $0\ls \vr(u_0)\ls 1$. Then,  for $\la>0$ small enough, by (3.18) and (3.20) we have
 $$\be{ll} u_0+\la(-V(u_0))&=(1-\la \vr(u_0)(1+\|u_0\|) u_0\\
 &+\la(1+\|u_0\|)\vr(u_0))\Big(u_0-\frac{1}{1+\|u_0\|} v(u_0)\Big)\in \Da.\en
 $$
 By  the arguments above and a well known theorem due to Brezis-Martin on flow invariant set,  we see that $\Da$ is a flow invariant set of (3.21).

 Now we show the following claim: \textit{for each $u_0\in A\cap J^{r+\va}\bk\mb K_r^2)_{3\de}$,  there exist an open neighborhood $U(x_0)$ of $x_0$ and  $\tau(u_0)>0$, such that $$\si(t, u)\in \mbox{int}(\Da)\cup \mbox{int} (\Db)\cup J^{r-\va}$$ for $t\gl t(u_0)$ and $u\in A\cap U(x_0)\cap J^{r+\va}$.}

  To prove the claim, we need to consider the following three cases:

 (1)\ \textit{$o(u_0)\cap\mbox{int}(\Da)\neq\ep$.}

  Assume without loss of generality that $u_1:=\si(t_1, u_0)\in \mbox{int}(\Da)$ for some $t_1\gl 0$. Take an open neighborhood  $U(u_1)$ of $u_1$ such that $U(u_1)\subset \mbox{int}(\Da)$. By the continuous dependence of ordinary differential
equations on initial data, there exists an open neighborhood $U(u_0)$ of $u_0$, and $\tau(u_0)>0$, such that $\si(t, u)\in U(u_1)\subset \mbox{int} (\Da)$ for all $u\in \big(A\cap J^{r+\va}\bk(\mb K_r^2)_{3\de}\big)\cap U(u_0)$ and $t\gl \tau(u_0)$.

 (2)\ \textit{ $o(u_0)\cap\mbox{int}(\Db)\neq\ep$.}

By the same way as the case (1) we can prove that, there exists an open neighborhood $U(u_0)$ of $u_0$, and $\tau(u_0)>0$, such that $\si(t, u_0)\in  \mbox{int}(\Db)$ for all $u\in \big(A\cap J^{r+\va}\bk(\mb K_r^2)_{3\de}\big)\cap U(u_0)$ and $t\gl \tau(u_0)$.

(3)\ \textit{$o(u_0)\cap\big(\mbox{int}(\Da)\cup \mbox{int}(\Da)\big)=\ep$, $o(u_0)\cap (\mb K^2_r)_{3\de}=\ep$.}

Let $h(t, u_0)=J(\si(t, u_0))$ for all $t\in [0,\infty)$. It is easy to see that $h(t, u_0)$ is locally Lipschitz continuous in $t\in [0,+\infty)$. According to Leburng's Mean Theorem  we have
$$\be{ll}\frac{\pa }{\pa s}h(s, u_0)&\ls \max\{\lge w^*, \frac{\pa}{\pa s} \si(s, u_0)\rge: w^*\in \pa J(\si(s, u_0))\}\ \mbox{a.e.}\\
  &=-\min\{ \lge w^*, V( \si(s, u_0))\rge: w^*\in \pa J(\si(s, u_0))\}\ \mbox{ a.e.}\\
 &\ls\left\{\be{ll}-\frac{b}{8}, &\mbox{if}\ \si(s, u_0)\in B(r,\bar\va, 2\de);\\
 0, &\mbox{otherwise}.\en\right.\en\no(3.23)$$
Consequently, we have for any $0\ls t_1<t_2<\infty$,
$$J(\si(t_1, u_0))-J(\si(t_2, u_0))=-\int^{t_2}_{t_1} \frac{\pa}{\pa s} h(s, u_0)ds\gl \frac{b}{8} (t_2-t_1) \ \mbox{if}\ \si(s, u_0)\in B(r, \va, 2\de)\ \mbox{and}\ s\in [t_1,t_2).\no(3.24)$$
So, if we let $\tau(u_0)>\frac{16 \va}{b}$, we have by (3.23),
$$J(\si(t, u_0))\ls J(u_0)-\frac{bt}{8}\ls r+\va-\frac{b}{8}<r-\va$$
for $t\gl \tau(u_0)$. Take an open neighborhood $U(\si(\tau(u_0),u_0))$ of $\si(\tau(u_0),u_0)$ such that $U(\si(\tau(u_0),u_0))\subset J^{r-\va}$. By the continuous dependence of ordinary differential
equations on initial data, there exists an open neighborhood $U(u_0)$ of $u_0$,  such that $\si(\tau(u_0), u)\in U(\si(\tau(u_0),u_0))$ for all $u\in \big(A\cap J^{r+\va}\bk(\mb K_r^2)_{3\de}\big)\cap U(u_0)$.  Since $h(t, u_0)$ is non-increasing in $t\in [0,+\infty)$,  we have  $\si(t, u)\in  J^{r-\va}$ for all $u\in \big(A\cap J^{r+\va}\bk (\mb K_r^2)_{3\de}\big)\cap U(u_0)$ and $t\gl \tau(u_0)$.

(4)\ \textit{$o(u_0)\cap\big(\mbox{int}(\Da)\cup \mbox{int}(\Da)\big)=\ep$, $o(u_0)\cap (\mb K^2_r)_{3\de}\neq\ep$.}

In this case we may take $[t_1,t_2]\subset [0,+\infty)$ such that $\si(t_1, u_0)\in \pa(K_r)_{3\de}$, $\si(t_2, u_0)\in \pa(K_r)_{2\de}$, and $\si(t, u_0)\in (K_r)_{3\de}\bk {(K_r)}_{2\de}$ for $t\in [t_1, t_2]$. By (3.21) we have
$$\be{ll}
 \|\si(t_2, u_0)-\si(t_1, u_0)\|&\ls \di\int^{t_2}_{t_1} \|V(\si(s, u_0)\| ds\\
 &\ls 2\di \int^{t_2}_{t_1}\big(1+\|\si(s, u_0)\|\big)ds\\
 &\ls 2\big(1+M_0\big)(t_2-t_1).
 \en\no (3.25)$$
So, we have
$$t_2-t_1\gl \frac{ \|\si(t_2, u_0)-\si(t_1, u_0)\|}{ 2\big(1+M_0\big)}\gl \frac{\de}{ 2\big(1+M_0\big)}.$$
 Now we show that there must exist $t_0\in [t_1,t_2]$ such that $\si(t_0,u_0)\in \{x\in X:J(x)<r-\va\}$. Assume by make contradiction that this is not the case, then $$\{\si(t,u_0):t\in [t_1,t_2]\}\cap \{x\in X: J(x)<r-\va\}=\ep.$$
 Then, by (3.24) and (3.25) we have
 $$\be{ll}J(\si(t_2, u_0))&\ls J(\si(t_1,u_0))-\frac{\ga}{4}(t_2- t_1)\\
 &\ls J(u_0)-\frac{\ga}{4}(t_2- t_1)\\
 &\ls r- \frac{\ga\de}{ 8\big(1+M_0\big)}\\
 &<r-\bar\va_0<r-\va,\en
$$
which is a contradiction. Thus, there must exist $t_0\in [t_1,t_2]$ such that $\si(t_0,u_0)\in \{x\in X:J(x)<r-\va\}$. Then, as the arguments in the case (3), we can show that there exists a open neighborhood of $U(u_0)$ of $u_0$, and $\tau(u_0)\gl 0$,  such that $\si(t, u)\subset \{x\in X:J(x)<r-\va\}$ for any $u\in U(u_0)$  and $t\gl \tau(u_0)$.

From the arguments above we see the claim  holds.  Obviously, $\mathscr U=\{U(u_0): u_0\in A\cap J^{r+\va}\bk\mb K_r^2)_{3\de}\}$ is an open cover of $A$. Since $A\cap J^{r+\va}\bk\mb K_r^2)_{3\de}$ is compact, there exist finite sets in $\mathscr U$, say $U(u_1), U(u_2),\cdots, U(u_m)$ such that $\big(A\cap J^{r+\va}\bk\mb K_r^2)_{3\de}\big)\subset\bigcup\limits_{i=1}^m U(u_i)$. Let $$T=\max\{\tau(u_1), \tau(u_2), \cdots, \tau(u_m)\}.$$ Then for each $u\in A\cap J^{r+\va}\bk\mb K_r^2)_{3\de}$, $\si(T, u)\in J^{-1}((-\infty, r-\va))\cup\Da\Db$. If we let $\eta(t, x)= \si(Tt,x)$ for $t\in [0,1]$ and $x\in X$. Then we can easily see the conclusions (1)$\sim$(5) hold.

Since $\Da\cup\Db$ is a flow invariant set of (3.21), we have $\eta(t, \Da\cup\Db)\subset \Da\cup\Db$. This means that the conclusion  (5) holds. The proof is complete.
\vskip 0.1in

{\bf Definition 3.1.}\ Let $Q, T$ be two closed subset of $X$, satisfies that (1) $\pa Q\cap T=\ep$; (2) for each $\ga\in C(Q, X)$ with $\ga|_{\pa Q}=I$, we have $\pa Q\cap T\neq \ep$. Then we call $(Q, T)$ is linking.
\vskip 0.1in

There have been many papers studied the critical points by using the linking structure. For Example, the paper [20] study the case where the functional $J$ being of $C^2$; while  [13] studied the case  where the functional $J$ is locally Lipschitz continuous.  Now we will use the linking structure and the deformation theorem to show the existence of sign-changing solutions. \vskip 0.1in

{\bf Theorem 3.1.}\  Suppose that (H$_1)$ and (H$_2)$ hold.   Let $W=\Da\cup\Db$. Let $(Q, T)$ is a linking. $T\subset X\bk W$ and $Q\subset X$ is compact, and for each $\ga\in\Gamma$, $\ga(Q)\cap S\neq \ep$, where $S=X\bk W$ and
$$\Gamma=\{\ga\in C(Q, X):\ga|_{\pa Q\cap S}=I, \ga(\pa Q\cap W)\subset W\}.$$
Assume that there exist $\beta>\al$ such that $$\sup\limits_{\pa Q\cap S}J(x)\ls \al<\beta\ls\inf\limits_T J(x)$$
and $\sup\limits_Q J(x)<\infty$. Let $r=\inf\limits_{\ga\in \Gamma}\sup\limits_{\ga(Q)\cap S} J(x)$. Then $r\gl \beta$ and $\mb K_r\cap S\neq\ep$.

{\bf Proof.}\ Since $Q\subset X$ is compact, $J$ is continuous on $Q$, then $J$ is bounded on $Q$. Since $I\in \Gamma$, we have
$$r\ls \sup\limits_{Q\cap S} J(x)\ls \sup\limits_QJ(x)<+\infty.$$
Now we show $ r\gl \beta$. Since for any $\ga\in \Gamma$, $\ga(Q)\cap T\neq\ep$ and $T\subset S$, then $\ga(Q)\cap S\cap T\neq \ep$. Thus, we have
$$\sup\limits_{\ga(Q)\cap S}J(x)\gl \sup\limits_{\ga(Q)\cap T\cap S} J(x)=\sup\limits_{\ga(Q)\cap T} J(x)\gl\inf\limits_{\ga(Q)\cap T} J(x)\gl \inf\limits_{T}J(x)\gl\beta.$$

Assume that $\mb K_r\cap \pa W=\ep$. Otherwise, the conclusion would hold. We claim that $\mb K_r\cap S\neq\ep$. Arguing by make contradiction that $\mb K_r\cap S=\ep$. Take $\va_0>0$ such that $r-2\va_0>\al$. Let $0<\va<\va_0$. Take $\ga\in\Gamma$ such that $\sup\limits_{\ga(Q)\cap S} J(x)\ls  r+\va$. Let $A=\ga(Q)$. By using the deformation theorem (Lemma 3.4),   there exists $\eta\in C([0,1]\ti X, X)$ such that the conclusions (1)$\sim$ (5) hold. In particular, we have $\eta(1, \ga(Q))\subset J^{r-\va}\cup W$.

Set $\ga_1=\eta(1, \ga(\cdot))$. Then  we have
$$\ga_1(\pa Q\cap W)=\eta(1, \ga(\pa Q\cap W))\subset \eta(1, W)\subset W.$$
For any $x\in \pa Q\cap S$  with $ J(x)\ls \al<r- 2\va$, we have $$\eta(1, \ga(x))=\eta(1, x)=x,$$
which implies that $\ga_1\in \Gamma$.

On the other hand, we have
 $$\sup\limits_{\ga_1(Q)\cap S} J(x)\ls \sup\limits_{(J^{r-\va}\cup W)\cap S} J(x)\ls\sup\limits_{J^{r-\va}} J(x)\ls r-\va.$$
which is a contradiction.  The proof is complete.
\vskip 0.1in

Now we apply  Theorem 3.1 to a special linking structure.
\vskip 0.1in

Suppose that $X=X_1\oplus X_1^\perp$, where $\infty>\mbox{dim}\ X_1\gl 1$. Take $y_1\in X_1^\perp$ with $\|y_1\|=1$. Let $ R>0$ and
$$Q=\{e= e_1+ty_1: e_1\in X_1,\|e_1\|\ls R, 0\ls t\ls R, \|e\|\ls R\}.$$
For  any $0<r< R$, let $T=\{ e\in X_1^\perp: \|e\|=r\}$. Then $(Q, T)$ is linking.
\vskip 0.1in

{ \bf Corollary 3.1.}\ Suppose that $J$ satisfies all conditions on Theorem 3.1, $(Q, T)$ is defined as above, $W\subset X_1^\perp\bk\{0\}$, and  for $\al<\beta$,
$$\max\limits_{\pa Q\cap S} J(x)\ls \al<\beta\ls\inf\limits_{T} J(x).\no(3.26)$$
 Let  $\Gamma, r$ be defined as Theorem 2.1. Then
 $$\mb K_r\cap (S\cup \pa W)\neq \ep.$$
 \vskip 0.1in

{\bf Remark 3.2.}\ Obviously, we can obtain the existence results of infinitely many critical points of Locally Lipschitz and even functionals by using the deformation lemma (Lemma 3.4). Here, for simplicity, we will not give the related  results. We think that readers of this paper can easily  obtain the related results using Lemma 3.4 by themselves.

\section{Sign Changing Solutions for  Differential Inclusion Problems}

\ \ \
Consider the following  differential  inclusions problems
$$\left\{\be{ll} -\De u \in \la \pa j(x, u) &\ \mbox{in}\ \Om,\\
u=0\ &\ \mbox{on}\ \pa\Om,\en\right.\no(4.1_\la)
$$
where $\Om $ is a bounded open domain in $\mb R^N$ with a $C^1$ boundary, $\la>0$ is a parameter, the reaction term $\pa j(x,s)$ is the generalized gradient of a non-smooth potential $s\mapsto j(x,s)$, which is subject to the following conditions.

{\bf (H$_j)$}\ $j:\Om\ti \mb R \mapsto \mb R$ is a Carath\'{e}odory  function and there exist constants $a_1>0$, $2<q<2^*$ such that

(i)\ $j(x,\cdot)$ is locally Lipschitz for almost every $x\in \Om$ and $j(\cdot, 0)=0$;

(ii)\ $|\xi|\ls a_1(1+|s|^{q-1})$ a.e. in $\Om$ and  for all $s\in \mb R$, $\xi\in\pa j(x,s)$;

(iii)\ For some $\mu>N(\frac{q}{2}-1)$ one has  $$\lim\inf\limits_{|x|\ri \infty}\frac{u\inf \pa j(x, z)-2j(x, z)}{|z|^\mu}> 0$$
uniformly  for almost $x\in \Om$;

(iv)\ $$\lim\sup_{z\ri 0}\frac{2j(x,z)}{z^2}= 0$$ uniformly with respect to $x\in \Om$;

(v)\ $zw\gl 0$ for each $w\in \pa j(x, z)$ a.e. $\Om$ and $z\in\mb R$.
\vskip 0.1in

{\bf Theorem 4.1.}\  Assume (H$_j)$ holds. Then (4.$1_\la$) has at least one  sign-changing solution $\bar u$ for $\la >0$.

\vskip 0.1in

Let $$\|u\|=\Big(\int_\Om|\na u|^2 dx\Big)^{1\over 2}, |u|_p=\Big(\int_\Om| u|^p dx\Big)^{1\over p}\ \mbox{and}\ |u|_\infty=\sup\limits_{x\in\Om}|u(x)|$$
be the standard norms of $H_0^1(\Om)$, $L^p(\Om)$, respectively  $C(\Om)$.
Let $X=H_0^1(\Om)$. We introduce the energy functional $J_\la: X\mapsto \mb R$ for $\la>0$ by
 $$J_\la (u)=\frac{1}{2}\|u\|^2-\la\vr(u),$$
 where  the functional $\vr: X\mapsto\mb R$ be defined by
$$\vr(u)=\int_\Om j(x, u(x))dx.$$
Let $P=\{u\in X: u(x)\gl 0\  \mbox{a.e.}\ \Om\}$. Let $\la_1,\la_2$ be the first and second eigenvalues of $-\Delta$ under the Dirichlet condition, and $u_1, u_2$ be the eigenfunctions corresponding to $\la_1$ and $\la_2$ respectively.
 \vskip 0.1in

{\bf Lemma 4.1.}\ Assume that (H$_j)$ holds. Then there exists $\mu_0>0$ such that $(I-\pa J_\la)\big({D}^\nu(\mu_0)\big)\subset D^\nu(\frac{1}{2}\mu_0)$
for $\nu\in\{+,-\}$.

{\bf Proof.}\ For each $u\in X$, denote $u^{\pm}=\max\{ \pm u, 0\}$. Then, for each $u\in X$ we have
$$\be{ll}|u^+|_2&=\min\limits_{w\in( -P)}|u-w|_2\\
&\ls \frac{1}{\sqrt{\la_1}}\min\limits_{w\in (-P)}\|u-w\|\\
&=\frac{1}{\sqrt{\la_1}}\mbox{dist}\ (u, -P).\en
$$
Similarly, we have for some constant $C_\mu>0$,
$$|u^+|_\mu\ls C_\mu \mbox{dist}\ (u, -P).$$

Obviously, $u\mapsto \frac{1}{2} \|u\|^2$ is a $C^1$-functional whose derivative is the  identity operator $I$. Aubin-Clarke's Theorem [4, Theorem 1.3.10] ensures  that the functional $\vr: X\ri\mb R$
is Lipschitz continuous on any bounded subset of $L^2(\Om)$ and its gradient is included in the set
$$N(u)=\{ w\in L^2(\Om): w(x)\in\pa j(x, u(x))\ \mbox{for almost every}\ x\in\Om\}.$$
Since $X$ continuously embedded in $L^2(\Om)$, the function $J_\la$ turns out to be locally Lipschitz on $X$. So, we have
$$ \pa J_\la (u)\subset \{u\}-\la N(u).\no(4.2)$$
Now, if $u\in X$  complies with $0\in \pa J_\la(u)$ then $$u=\la w\ \mbox{in}\ X^*$$
for some $w\in N(u)$.

By (H$_j)$ (i) and (iv), for any $u\in X$, $z(x)\in \pa j(x, u(x))$, we have
$$( z(x), u(x))\ls \la\va |u(x)|^2+\la C_\va|u(x)|^{q}.$$
Now we check the Schauder invariance condition. Take $\mathscr F: X^*\ri X$ as the identity operator $I:X\ri X$. By (4.2) we have
$$(I-\pa J_\la)(u)\subset \la N(u).$$
Let $w_\la\in\la N(u)$ with $w_\la=\la w$ for some $w\in N(u)$. So, if $\va=\va(\la)>0$ small enough,
$$\be{ll}\mbox{dist} (w_\la, -P)\|w_\la^+\|&\ls \|w_\la^+\|^2=(w_\la, w_\la^+)\\
&\ls\di\int_\Om \big(\la\va |u^+|+\la C_\va|u^+|^{q-1}\big)|w^+_\la| dx\\
&\ls \Big(\frac{1}{3} \mbox{dist}\ (u, -P)+C\big(\mbox{dist}\ (u, -P)\big)^{q-1}\Big)\|w^+_\la\|,\en
$$
where $C>0$ is a constant, here we have employed some method of [25,p.71]. Thus, we have
$$\mbox{dist}\ \big((I-\pa J_\la)(u), -P\big)\ls \frac{1}{3} \mbox{dist}\ (u, -P)+C\big(\mbox{dist}\ (u, -P)\big)^{q-1}.$$
Hence, there exists $\mu_1>0$ such that for all $0<\bar\mu<\mu_1$
$$(I-\pa J_\la)\big({D}^-(\bar\mu)\big)\subset D^-(\frac{1}{2}\bar\mu).$$

Similarly, we can show that  there exists a $\mu_2>0$ such that for all $0<\bar\mu<\mu_2$,
$$(I-\pa J_\la)\big({D}^+(\bar\mu)\big)\subset D^+(\frac{1}{2}\bar\mu).$$
Let $\mu_0=\min\{\mu_1,\mu_2\}$. Then we have
$$(I-\pa J_\la)\big({D}^\nu(\mu_0)\big)\subset D^\nu(\frac{1}{2}\mu_0)$$
for $\nu\in\{+,-\}$. So, the  condition  (H$_1)$ holds. The proof is complete.

\vskip 0.1in

Now we make a decomposition for the real Banach space $X: X=\mbox{span}\ \{u_1\}\oplus V$, where $V=\{ v\in X: \int_\Om v u_1=0\}$.
\vskip 0.1in
\vskip 0.1in

{\bf  Lemma 4.2.}\ Assume (H$_j)$ holds. Then $J_\la$ satisfies the condition (H$_2)$.

{\bf Proof. }\ The proof is similar to Lemma 5.1 in [11] and Proposition 7 in [6]. For reader's convenience we give the details of the process. We divide the proof into three steps.

{\bf Step 1.}\  $J_\la$ satisfies the condition  $(PS)_X$.

Let $\{ x_n\}\subset X$  be such that $|J_\la(x_n)|\ls M_1$ for some $M_1>0$ and $(1+\|x_n\|)m(x_n)\ri 0$ as $n\ri \infty$.
Take $x_n^*\in \pa J_\la(x_n)$ such that $m(x_n)=\|x_n^*\|_{*}$ for $n\gl 1$. Then we have $x_n^*=A (x_n)- u_n$ with $A: X\ri X^*$ being defined by
$$( A(u), y)=\int_\Om \na u\cdot \na y dx \ \ \mbox{for}\ u,y\in X,$$
and $u_n\in L^2(\Om)$, $u_n(x)\in \pa \vr(x, x_n(x))$. Since $(1+\|x_n\|)m(x_n)\ri 0$ as $n\ri \infty$, we can say that $|( x_n^*, x_n)|\ls \frac{1}{n}$. So,
$$-\|x_n\|^2+\la(u_n,x_n)\ls \frac{1}{n}.\no(4.3)$$

Similarly, since $|2 J_\la (x_n)|\ls 2M_1$ for all $n\gl 1$, we have
$$\|x_n\|^2-2\la\int_\Om j(x, x_n(x))dx\ls 2M_1.\no(4.4)$$
It follows from (4.3) and (4.4) that
$$\int_\Om\big(\la u_n(x)\cdot x_n(x)-2 j(x, x_n(x))\big) dx\ls \frac{1}{n}+2M_1. \no(4.5)$$

By the condition (H$_j)$ (iii), there exists $M_2>0$ such that for a.e. $x\in \Om$ and $|z|\ls M_2$, and all $u\in \pa j(x, z)$, we have
$$\frac{\beta}{2}|z|^\mu< uz- 2 j(x, z).$$
By using (H$_j)$ (ii),  for a.e. $x\in \Om$ and $|z|\ls M_2$, and all $u\in \pa j(x, z)$, we have
 $$\frac{\beta}{2}|z|^\mu- a_2(x)< uz- 2 j(x, z), a_2\in L^\infty(\Om), a_2(x)\gl 0 \ \mbox{a.e. on }\ \pa\Om.\no(4.6)$$
By (4.4) and (4.5) we have
 $$\frac{\beta}{2}|x_n|^\mu_\mu\ls\int_\Om( u_n(x)x_n(x)- 2 j(x, x_n(x)))dx+| a_2(x)|_1<\frac{1}{n}+2M_1+|a_2|_\infty.$$
 This implies  that $\{|x_n|_\mu\}$ is bounded.

 Take $s>1$ such that $2<s<\min\big\{ 2^*, \frac{2(N+\mu)}{N}\big\}$,
 and for a.e. $x\in \Om$, $z\in \mb R$,
 $$j(x, z)\ls c_2+c_3 |z|^s, c_2, c_3>0.$$
 Take $$\theta=\left\{\be{ll}\frac{2^*(s-\mu)}{s(2^*-\mu)}, &{N>2},\\
 1-\frac{\mu}{s},& N=2.\en\right.
 $$
 Then we have $0<\theta<1$, and $$\frac{1}{s}=\frac{1-\theta}{\mu}+\frac{\theta}{2^*}.$$

 Thus, by using the interpolation inequality and the Sobolev embedding Theorem, we obtain
 $$|x_n|_s\ls |x_n|_\mu^{1-\theta}|x_n|_{2^*}^\theta\ls c_4 \|x_n\|^\theta \ \mbox{for some }\ c_4>0.$$
 Since $J_\la(x_n)\ls M_1$ and $L^s(\Om)\hookrightarrow L^2(\Om)$, by using the well known Young's inequality we have
 $$\be{ll}\frac{1}{2}\|x_n\|^2&\ls c_2|\Om|+c_3|x_n|_s^s+M_1\\
&\ls c_5+c_6|x_n|_s^s\ls c_7+c_8\|x_n\|^{\theta s},\en
$$
where $c_5,\cdots, c_8$ are positive constants.

If $N>2$, then $Ns <2(N+\mu)$,  and so we have
$$\theta s=\frac{2^*(s-\mu)}{2^*-\mu}=\frac{2N}{N-2}\cdot\frac{(s-\mu)(N-2)}{2N-N\mu+2\mu}<\frac{2N}{N-2}\cdot\frac{(s-\mu)(N-2)}{Ns-N\mu}=2.$$
If $N=2$, then $s<\min\{ 2^*, \frac{2(N+\mu)}{N}\}= 2+\mu$. So,
$$\theta s=\big(1-\frac{\mu}{s}\big)s=s-\mu <2.$$
Thus,  we have $\theta s<2$ for $N\gl 2$. Thus, $\{\|x_n\|\}$ is bounded.

Assume that $x_n\rightharpoonup x$ in $X$ and $x_n\ri x$ in $L^2(\Om)$. By usual method we can show that $x_n\ri x$ in $X$. So, $J_\la$ satisfies the condition on $(PS)_X$.

{\bf Step 2.}\  $J_\la$ satisfies the condition  $(PS)_{\Da, D^+(\frac{3}{4}\mu_0)}$.

 Let $\{ x_n\}\subset D^+(\frac{3}{4}\mu_0)$  be such that $|J_\la(x_n)|\ls M_2$ for some $M_2>0$,  and $$(1+\|x_n\|)m_{\Da}(x_n)\ri 0\ \mbox{as}\ n\ri \infty.$$
Let $h_{\Da}: X^*_w\ti X\mapsto \mb R$ be defined in the manner in [6]  by
$$h_{\Da}(x^*,x)=\sup\big\{\lge x^*, x-y\rge: y\in \Da, \|y-x\|<1\big\}.$$
 Then, we can find  $x_n^*\in \pa J_\la(x_n)$ such that $m_{\Da}(x_n)= h_{\Da}(x_n^*, C_1(x_n))$, where $C_1(x_n)=(x_n-\Da)\cap B(0,1)$. That is
$$( x_n^*, x_n-y)\ls h_{\Da}(x_n^*, C_1(x_n))\ \mbox{for all}\ y\in\Da, \|x_n-y\|<1.\no(4.7)$$
Obviously, $x_n^*=A (x_n)- \la u_n$ with $A: X\ri X^*$ being defined as Step 1.  Let us denote by $\pi_P: X\mapsto P$ the orthogonal projection on the cone $P$ so that $T_P(u)=u-\pi(u)$ realizes the distance of $u$ from $P$. Let $y_n=x_n+\frac{\mu_0}{8\|x_n\|} x_n$ for $n\gl 1$. Since $\pi(x_n)\in P$, we have
$$\be{ll}\mbox{dist}\ (y_n, P)&=\|T_Py_n\|\ls \|y_n-\pi(x_n)\|\\
&\ls \|y_n-x_n\|+\|x_n-\pi(x_n)\|=\|y_n-x_n\|+\|T_px_n\|\\
&=\|y_n-x_n\|+\mbox{dist}\ (x_n, P)\ls \frac{\mu_0}{8}+\frac{3\mu_0}{4}\ls \mu_0.\en$$ Thus, $y_n\in \Da$, $\|x_n-y_n\|<1$. And so,  by (4.7) we have
$$\be{ll} (1+\|x_n\|)(x_n^*, x_n- y_n)&=-\frac{\mu_0(1+\|x_n\|)}{8\|x_n\|}(x_n^*,x_n)\\
&=-\frac{\mu_0(1+\|x_n\|)}{8\|x_n\|}\Big((A(x_n),x_n)-\la\di\int_\Om u_n x_n dx\Big)\\
&\ls (1+\|x_n\|)h_{\Da}(x_n^*, C_1(x_n))=:\va_n\ \mbox{with}\ \va_n\downarrow 0.\en
$$
So we obtain
$$-\|x_n\|^2+\la \int_\Om u_n x_n\ls \frac{8\|x_n\|}{\mu_0(1+\|x_n\|)}\va_n=\va'_n, \va_n'\downarrow 0.$$
Then, by using the arguments from (4.3) to the end of the step 1 we can show that $\{\|x_n\|\}$ is bounded, and $\{x_n\}$ has a convergent subsequence.

{\bf Step 3.}\  $J_\la$ satisfies the condition $(PS)_{\Db, D^-(\frac{3}{4}\mu_0)}$.

In a similar way as Step 2, we can show the claim above hold.

{\bf Step 4.}\  $J_\la$ satisfies the condition $(PS)_{\Dc}$.

Obviously,  $\Dc$ is bounded in $X$. So, by using the Sobolev inequality,  we  see that $\Dc$ is bounded in  $L^\mu(\Om)$ for $\mu\in (1, 2^*)$. Then, by using the arguments above we can show that $J_\la$ satisfies the condition $(PS)_{\Dc}$.

The proof is complete.
\vskip 0.1in

{\bf Lemma 4.3.}\ Assume that (H$_j)$ holds. Then there exists $\de>0$ such that $J_\la(v)> 0$ for all $v\in V$ with $0<\|v\|\ls \de$.

{\bf Proof.}\ For each $\va>0$,  there exists $C_1(\va)>0$ such that
$$j(x, z)\ls \frac{\va}{2}|z|^2+C_1|z|^q.$$
So, for any $v\in V$ we have
$$\be{ll}
J_\la(v)&=\frac{1}{2}\|v\|^2-\la\di\int_\Om j(x, v(x))dx\\
&\gl \frac{1}{2}\|v\|^2-\la\di\int_\Om\big(\frac{\va}{2}|v(x)|^2+C_1|v(x)|^q\big) dx\\
& \gl \frac{1}{2}\big(1-\frac{\la}{\la_V}\va)\|v\|^2-C_2\|v(x)\|^q,\en$$
where $C_2>0$ is a constant. Take $\va>0$ small such that $\frac{\la}{\la_V}\va<1$. Then, we have $J_\la(v)> 0$ for all $v\in V$ with $0<\|v\|\ls \de$. The proof is complete.
\vskip 0.1in

{\bf Lemma 4.4.}\ Assume that (H$_j)$ holds. Then there exists $R>0$ such that $J_\la(u)< 0$ for all $u\in E_2$ with $\|u\|\gl R$, where $E_2=\mbox{span}\ \{\vr_1,\vr_2\}$.

{\bf Proof.}\ For each $u\in E_2$, we have $\|u\|^2<\la_2|u|_2^2$. By the condition (H$_j) (iii)$ we have
$$j(x, u)\gl a_0|u|^\mu-C_3$$
for some $a_0, C_3>0$. So, we have
$$\be{ll} J_\la (u)&=\frac{1}{2}\|u\|^2-\la\di\int_\Om j(x, u)dx\\
&\ls \frac{1}{2}\|u\|^2-\la a_0|u|^\mu_\mu+C_4\\
&\ls \frac{1}{2}\|u\|^2-\la C_\mu\|u\|^\mu+C_4,\en
$$
where $C_4>0$ is a constant. Then we have
$$\lim\limits_{ u\in E_2, \|u\|\ri\infty} J_\la (u)=-\infty,$$
and so the conclusion holds. The proof is complete.
\vskip 0.1in

{\bf The Proof of Theorem 4.1.}\ We will apply Corollary 3.1 to show Theorem 4.1.  It follows from  Lemma 4.1 and 4.2 that the conditions (H$_1)$ and (H$_2)$ hold.  Take $R>0$ large enough such as that in Lemma 4.4. Let
$$Q=\{e= \vr_1+t\vr_2\in E_2: e_1\in X_1,\|e_1\|\ls R, 0\ls t\ls R, \|e\|\ls R\},$$
and $T=\{x\in V:\|x\|=\de\}$, where $\de>0$ be defined as that in Lemma 4.3. It follows from Lemma 4.3 and 4.4 that (3.26) hold. Hence, by Corollary 3.1, we see that Theorem 4.1 hold. The proof is complete.
\vskip 0.1in

{\bf Remark 4.1.}\ There have been some papers studied the existence for sign-changing solutions of differential inclusion problems; see [14,15] and the references therein. For example, by combing variational methods with truncation techniques the paper [14] obtained the existence of positive, negative and nodal solutions to differential inclusion problems with a parameter. Here, our  method is different to  that in [14,15].

\end{document}